\documentclass[11pt,oneside]{article}
\usepackage{amsfonts,amssymb,amsmath,verbatim,authblk}
\usepackage{hyperref}

\newcommand{\beq}{\begin{equation}\ }
\newcommand{\eeq}{\end{equation}\ }

\begin{document}

\author[1,2]{Demetris T. Christopoulos}
\affil[1]{\small{National and Kapodistrian University of Athens, Department of Economics}}
\affil[2]{\tt{dchristop@econ.uoa.gr}, \tt{dem.christop@gmail.com} }

\title{{\itshape Linear Regression without computing pseudo-inverse matrix}}

\maketitle

\begin{abstract}
We are presenting a method of linear regression based on Gram-Schmidt orthogonal projection that does not compute a pseudo-inverse matrix. This is useful when we want to make several regressions with random data vectors for simulation purposes.\medskip
\end{abstract}


\section{The traditional OLS regression solution}
We have the data
\beq
y_i=\beta_0+\beta_1\,x_i+\epsilon_i,\,i=1,2,\ldots,n,\,\,\epsilon_i\sim iid\left(0,\sigma^2\right)
\eeq
In matrix form we follow \cite{hay-00} and can write
\beq
y=\mathbf{X}\beta+\epsilon,\,\,\epsilon\sim iid\left(0,\sigma^2\,I_2\right)
\eeq
with
\beq
{\underset{(2\times 1)}{\mathbf{x}_{i}}}=\begin{bmatrix}1\\ x_{i}\end{bmatrix}
,\,{\underset{(2\times 1)}{\beta}}=\begin{bmatrix}\beta_{0}\\ \beta_{1}\end{bmatrix},\,	
\eeq
and
\beq
\underset{(n\times 1)}{y}=\begin{bmatrix}y_1\\ y_2 \\ \vdots \\ y_n\end{bmatrix},\,
\underset{(n\times 1)}{\epsilon}=\begin{bmatrix}\epsilon_1\\ \epsilon_2 \\ \vdots \\ \epsilon_n  \end{bmatrix},\,
\underset{(n\times 2)}{\mathbf{X}}=\begin{bmatrix}{\mathbf{x}}_{1}^{'}\\ {\mathbf{x}}_{2}^{'} \\ \vdots \\ {\mathbf{x}}_{n}^{'} \end{bmatrix}=
\begin{bmatrix}1&x_{1} \\ 1&x_{2} \\ \vdots \\  1&x_{n}\end{bmatrix}	
\eeq
The traditional Ordinary Least Squares solution is given by 
\beq
\hat{\beta}=\left(\mathbf{X}^{'}\mathbf{X}\right)^{-1}\mathbf{X}^{'}y=\mathbf{X}^{+}y
\eeq
where the Moon-Penrose pseudo-inverse matrix is defined, after \cite{moo-20} and \cite{pen-55}, as
\beq
\mathbf{X}^{+}=\left(\mathbf{X}^{'}\mathbf{X}\right)^{-1}\mathbf{X}^{'}
\eeq
The solution $\hat{\beta}$ can be written in terms of sums as 
\beq
\label{eq:sum}
\begin{bmatrix}\hat{\beta_{0}}\\ \\ \hat{\beta_{1}}\end{bmatrix}= \\
\begin{bmatrix}
\bar{y}-\hat{\beta_{1}}\bar{x}\\ \\
\frac{\sum_{i=1}^{n}{\left(x_{i}-\bar{x}\right)\,\left(y_{i}-\bar{y}\right)}}{\sum_{i=1}^{n}{\left(x_{i}-\bar{x}\right)^2}}
\end{bmatrix}=
\begin{bmatrix} 
\bar{y}-\frac{\sum_{i=1}^{n}{\left(x_{i}-\bar{x}\right)\,\left(y_{i}-\bar{y}\right)}}{\sum_{i=1}^{n}{\left(x_{i}-\bar{x}\right)^2}}\,\bar{x}\, \\ \\
\frac{\sum_{i=1}^{n}{\left(x_{i}-\bar{x}\right)\,\left(y_{i}-\bar{y}\right)}}{\sum_{i=1}^{n}{\left(x_{i}-\bar{x}\right)^2}}
\end{bmatrix}
\eeq

We want to avoid the computation of $\left(\mathbf{X}^{'}\mathbf{X}\right)^{-1}$ matrix. For this reason we recall from the geometry of least squares, \cite{joh-97}, that OLS linear regression is merely an orthogonal projection of data vector y in the column space of $\mathbf{X}$. 
In order to do such a projection we first perform a Gram-Schmidt orthonormalization process.
Our first vector is the first column of ones of the matrix $\mathbf{X}$
\beq
v_1=X_1=\begin{bmatrix}1\\1\\ \vdots \\1 \end{bmatrix}
\eeq
Our second vector is 
\beq
v_2=X_2-\underset{v_1}{proj}{X_2}=X_2-\frac{X_2^{'}v_1}{v_1^{'}v_1}\,v_1=
\begin{bmatrix}x_{1}\\x_{2}\\ \vdots \\x_{n} \end{bmatrix}-\frac{\sum_{i=1}^{n}{x_{i}}}{n}\,\begin{bmatrix}1\\1\\ \vdots \\1 \end{bmatrix}=
\begin{bmatrix}x_{1}-\bar{x}\\x_{2}-\bar{x}\\ \vdots \\x_{n}-\bar{x} \end{bmatrix}
\eeq
The two vectors $v_i,i=1,2$ are orthogonal. We normalize them by dividing with their norm and we have the vectors
\beq
\xi_1=\frac{1}{\sqrt{n}}\,\begin{bmatrix}1\\1\\ \vdots \\1 \end{bmatrix}
\eeq
\beq
\xi_2=\frac{1}{\sqrt{\sum_{i=1}^{n}{\left(x_{i}-\bar{x}\right)^2}}}\,\begin{bmatrix}x_{1}-\bar{x}\\x_{2}-\bar{x}\\ \vdots \\x_{n}-\bar{x} \end{bmatrix}
\eeq
Now the orthogonal projection of y in the vector space defined by $\xi_1,\xi_2$ is
\beq
\hat{y}=\left(y\cdot\xi_1\right)\xi_1+\left(y\cdot\xi_2\right)\xi_2
\eeq
where $u\cdot v$ denotes the inner product of vectors $u,v$. We find that
\beq
\hat{y}=\frac{1}{\sqrt{n}}\,\sum_{i=1}^{n}{y_i}\,\begin{bmatrix}1\\1\\ \vdots \\1 \end{bmatrix}+
\frac{\sum_{i=1}^{n}{\left(x_{i}-\bar{x}\right)\,y_i} }{\sqrt{\sum_{i=1}^{n}{\left(x_{i}-\bar{x}\right)^2}}}\,    
\frac{1}{\sqrt{\sum_{i=1}^{n}{\left(x_{i}-\bar{x}\right)^2}}}\,\begin{bmatrix}x_{1}-\bar{x}\\x_{2}-\bar{x}\\ \vdots \\x_{n}-\bar{x} \end{bmatrix}
\eeq
or
\beq
\hat{y}=\bar{y}\,\begin{bmatrix}1\\1\\ \vdots \\1 \end{bmatrix}+
\frac{\sum_{i=1}^{n}{\left(x_{i}-\bar{x}\right)\,y_i} }{\sum_{i=1}^{n}{\left(x_{i}-\bar{x}\right)^2}}\,    
\begin{bmatrix}x_{1}-\bar{x}\\x_{2}-\bar{x}\\ \vdots \\x_{n}-\bar{x} \end{bmatrix}
\eeq
By re-arranging the terms and by concerning the column vectors of $\mathbf{X}$ we obtain
\beq
\hat{y}=\left(\bar{y}-\frac{\sum_{i=1}^{n}{\left(x_{i}-\bar{x}\right)\,y_i}}{\sum_{i=1}^{n}{\left(x_{i}-\bar{x}\right)^2}}\,\,\bar{x}\right)\,
\begin{bmatrix}1\\1\\ \vdots \\1 \end{bmatrix}+
\frac{\sum_{i=1}^{n}{\left(x_{i}-\bar{x}\right)\,y_i} }{\sum_{i=1}^{n}{\left(x_{i}-\bar{x}\right)^2}}\,    
\begin{bmatrix}x_{1}\\x_{2}\\ \vdots \\x_{n}\end{bmatrix}
\eeq
or by recalling that
$$
\sum_{i=1}^{n}{\left(x_{i}-\bar{x}\right)}=0
$$
we finally have that the orthogonal projection is
\beq
\hat{y}=\left(\bar{y}-\frac{\sum_{i=1}^{n}{\left(x_{i}-\bar{x}\right)\,\left(y_{i}-\bar{y}\right)}}{\sum_{i=1}^{n}{\left(x_{i}-\bar{x}\right)^2}}\,\,\bar{x}\right)\,
\begin{bmatrix}1\\1\\ \vdots \\1 \end{bmatrix}+
\frac{\sum_{i=1}^{n}{\left(x_{i}-\bar{x}\right)\,\left(y_{i}-\bar{y}\right)}}{\sum_{i=1}^{n}{\left(x_{i}-\bar{x}\right)^2}}\,
\begin{bmatrix}x_{1}\\x_{2}\\ \vdots \\x_{n}\end{bmatrix}
\eeq
which gives directly the results of Eq. \ref{eq:sum} for $\hat{\beta_{0}},\hat{\beta_{1}}$.
So the two methods give the same results for the OLS solution $\hat{y}=\hat{\beta_{0}}+\hat{\beta_{1}}\,\mathbf{x}$.\\
Obviously the above technique can be generalized for every $n\times{k}$ matrix $\mathbf{X}$ and gives the projected vector $\hat{y}=\mathbf{X}\,\hat{\beta}$, although it is not easy to compute the components of $\hat{\beta}$ for the case $k>2$. For simulation purposes we only need the projected vector $\hat{y}$ and not really $\hat{\beta}$, thus the method is useful for those cases.


\begin{thebibliography}{3}
\bibitem{hay-00}{\sc F. Hayashi},
{\em Econometrics}, Princeton University Press. 2000

\bibitem{moo-20}{\sc E.~ H. Moore}, On the reciprocal of the general algebraic matrix.
{\em Bulletin of the American Mathematical Society} \textbf{26} (9): 394–395, 1920

\bibitem{pen-55}{\sc R. Penrose}, A generalized inverse for matrices.
{\em Proceedings of the Cambridge Philosophical Society} \textbf{51}: 406–413, 1955


\bibitem{joh-97}{\sc J. Johnston ans J. DiNardo},
{\em Econometric Methods}, McGraw-Hill, 4th ed. 1997


\end{thebibliography}
\end{document}